\documentclass[12pt]{article}
\usepackage{amsmath,amsfonts,amssymb}

\def\qed{\nopagebreak\hfill{\rule{4pt}{7pt}}}
\def\proof{\noindent {\it{Proof.} \hskip 2pt}}

\def\Z{\mathbb Z}

\newtheorem{theo}{Theorem}[section]
\newtheorem{deff}[theo]{Definition}
\newtheorem{lemm}[theo]{Lemma}
\newtheorem{prop}[theo]{Proposition}
\newtheorem{coro}[theo]{Corollary}

\begin{document}

\begin{center}
{\large \bf The Zrank Conjecture and Restricted Cauchy Matrices}
\end{center}

\begin{center}
Guo-Guang Yan$^{1}$, Arthur L. B. Yang$^{2}$ and Joan J.
Zhou$^{3}$\\[6pt]
$^{1}$Information School, Zhongnan University of Economics and Law\\
Wuhan 430060, P. R. China\\
Email: $^{1}${\tt guogyan@eyou.com}\\[6pt]
$^{1,2,3}$Center for Combinatorics, LPMC\\
Nankai University, Tianjin 300071, P. R. China\\
Email: $^{2}${\tt yang@nankai.edu.cn}, $^{3}${\tt jinjinzhou@hotmail.com}\\
\end{center}

\centerline{March 28, 2005}

\vspace{0.3cm} \noindent{\bf Abstract.} The rank of a skew
partition $\lambda/\mu$, denoted ${\rm rank}(\lambda/\mu)$, is the
smallest number $r$ such that $\lambda/\mu$ is a disjoint union of
$r$ border strips. Let $s_{\lambda/\mu}(1^t)$ denote the skew
Schur function $s_{\lambda/\mu}$ evaluated at
$x_1=\cdots=x_t=1,\,x_i=0$ for $i>t$. The zrank of $\lambda/\mu$,
denoted ${\rm zrank}(\lambda/\mu)$, is the exponent of the largest
power of $t$ dividing $s_{\lambda/\mu}(1^t)$. Stanley conjectured
that ${\rm rank}(\lambda/\mu)={\rm zrank}(\lambda/\mu)$. We show
the equivalence between the validity of the zrank conjecture and
the nonsingularity of restricted Cauchy matrices. In support of
Stanley's conjecture we give affirmative answers for some special
cases.


\noindent {\bf MSC2000 Subject Classification:} 05E10, 15A15.

\noindent {\bf Suggested Running Title:}  The Zrank Conjecture

\noindent {\bf Corresponding Author:} Arthur L. B. Yang,
yang@nankai.edu.cn

\section{Introduction}

Let $\lambda=(\lambda_1,\,\lambda_2,\ldots)$ be a partition of an
integer $n$, i.e., $\lambda_1\geq \lambda_2\geq \cdots \geq 0$ and
$\lambda_1+\lambda_2+\cdots=n$. The number of positive parts of
$\lambda$ is called the length of $\lambda$, denoted
$\ell(\lambda)$. The \emph{Young diagram} of $\lambda$ may be
defined as the set of points $(i,j)\in \Z^2$ such that $1\leq
j\leq\lambda_i$ and $1\leq i\leq \ell(\lambda)$. A Young diagram
can also be represented in the plane by an array of squares
justified from the top and left corner with $\ell(\lambda)$ rows
and $\lambda_i$ squares in row $i$. A square $(i,j)$ in the
diagram is the square in row $i$ from the top and column $j$ from
the left. The content of $(i,j)$, denoted $\tau((i,j))$, is given
by  $j-i$. The \emph{rank} of $\lambda$, denoted ${\rm
rank}(\lambda)$, is the length of the main diagonal of the diagram
of $\lambda$. Given two partitions $\lambda$ and $\mu$, we say
that $\mu\subseteq\lambda$ if $\mu_i\leq \lambda_i$ for all $i$.
If $\mu\subseteq\lambda$, we define a \emph{skew partition}
$\lambda/\mu$, whose Young diagram is obtained from the Young
diagram of $\lambda$ by peeling off the Young diagram of $\mu$
from the upper left corner.

We assume that the reader is familiar with the notation and
terminology on symmetric functions  in \cite{S1}. In connection
with tensor products of Yangian modules, Nazarov and Tarasov
\cite{NT} give a generalization of a rank to a skew partition
$\lambda/\mu$. Recently Stanley developed a general theory of
minimal border strip decompositions and gave several simple
equivalent characterizations of ${\rm rank}(\lambda/\mu)$ in
\cite{S2}. One of the characterizations of the rank of a skew
partition $\lambda/\mu$ says that  ${\rm rank}(\lambda/\mu)$ is
the smallest integer $r$ such that the Young diagram of
$\lambda/\mu$ is the disjoint union of $r$ border strips. Let
$s_{\lambda/\mu}(1^t)$ denote the skew Schur function
$s_{\lambda/\mu}$ evaluated at $x_1=\cdots=x_t=1,\,x_i=0$ for
$i>t$. The \emph{zrank} of $\lambda/\mu$, denoted ${\rm
zrank}(\lambda/\mu)$, is the largest power of $t$ dividing the
polynomial $s_{\lambda/\mu}(1^t)$.  Stanley conjectured that the
equality ${\rm rank}(\lambda/\mu)={\rm zrank}(\lambda/\mu)$ always
holds, which we call the \emph{zrank conjecture}.

In his combinatorial approach to the zrank conjecture in
\cite{S2}, Stanley defined the snake sequence and the interval
sets for a skew partition $\lambda/\mu$. In Section 2 for each
interval set $\mathcal{I}$ of $\lambda/\mu$ we define an interval
permutation $\sigma_{\mathcal{I}}$. Let ${\rm cr}(\mathcal{I})$ be
the number of crossings of $\mathcal{I}$, and let
$\rm{inv}(\sigma_{\mathcal{I}})$ be the number of inversions of
$\sigma_{\mathcal{I}}$. We show that ${\rm cr}(\mathcal{I})$ and
$\rm{inv}(\sigma_{\mathcal{I}})$ have the same parity.

Stanley generalized the code of a partition to the code of a skew
partition, and obtained a two-line binary sequence in \cite{S2}.
This sequence is called the \emph{partition sequence} by
Bessenrodt \cite{B1, B2}. Given a minimal border strip
decomposition $\mathbf{D}$ of $\lambda/\mu$, let $P_{\mathbf{D}}$
be the set of the contents of the lower left-hand squares of the
border strips in $\mathbf{D}$, and let $Q_{\mathbf{D}}$ be the set
of the contents of the upper right-hand squares. Using the
partition sequence, we show that $P_{\mathbf{D}}$ and
$Q_{\mathbf{D}}$ are uniquely determined by the shape of the skew
partition $\lambda/\mu$ in Section 3, i.e., these two sets are
independent of the minimal border strip decomposition
$\mathbf{D}$. For a given skew partition, we find a connection
between the values of these two sets and the paired integers of
the interval set.

Outside decompositions are introduced by Hamel and Goulden
\cite{HG} and are used to give a unified approach to the
determinantal expressions for the skew Schur funtions including
the Jacobi-Trudi determinant, its dual, the Giambelli determinant
and the ribbon determinant. For any outside decomposition, Hamel
and Goulden derive a determinantal formula with ribbon Schur
functions as entries. Their proof is based on a lattice path
construction and the Gessel-Viennot methodology \cite{GV1, GV2}.
In Section 4 we employ the determinantal formula in the case of
the greedy border strip decomposition and give the evaluation of
$(t^{-{\rm rank}(\lambda/\mu)}s_{\lambda/\mu}(1^t))_{t=0}$. As a
consequence we obtain the combinatorial description of $(t^{-{\rm
rank}(\lambda/\mu)}s_{\lambda/\mu}(1^t))_{t=0}$ in terms of the
interval sets of $\lambda/\mu$ given by Stanley \cite[Eq.
(30)]{S2}.

Based on the above results, we give an equivalent characterization
of the zrank conjecture. Given two positive integer sequences, we
define a \emph{restricted Cauchy matrix} corresponding to these
two sequences. The main objective of this paper is to show that
 the zrank conjecture holds for any skew partition  if and only if
all the restricted Cauchy matrices are nonsingular.  We present a
constructive proof for this equivalence in Section 5. Using some
fundamental properties of determinants, we confirm the
nonsingularity of the restricted Cauchy matrices for several
special classes of skew partitions.

\section{Snake sequences and interval sets}

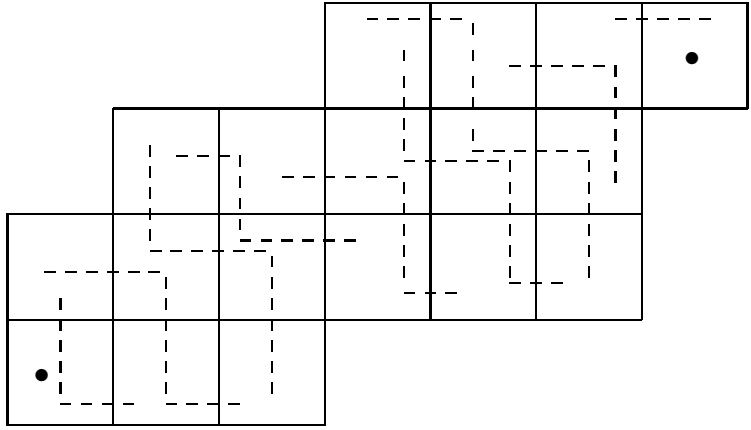
\begin{figure}[h,t]
\setlength{\unitlength}{2pt}
\begin{center}
\begin{picture}(160,100)
\put(80,85){\line(1,0){80}} \put(40,65){\line(1,0){120}}
\put(20,45){\line(1,0){120}} \put(20,25){\line(1,0){120}}
\put(20,5){\line(1,0){60}}

 \put(20,5){\line(0,1){40}}
\put(40,5){\line(0,1){60}} \put(60,5){\line(0,1){60}}
\put(80,5){\line(0,1){80}} \put(100,25){\line(0,1){60}}
\put(120,25){\line(0,1){60}} \put(140,25){\line(0,1){60}}
\put(160,65){\line(0,1){20}}

\put(25,13){$\bullet$} 

\put(30,9){$\line(1,0){2}$ $\line(1,0){2}$ $\line(1,0){2}$
$\line(1,0){2}$}

\put(30,11){$\line(0,1){2}$} \put(30,15){$\line(0,1){2}$}
\put(30,19){$\line(0,1){2}$} \put(30,23){$\line(0,1){2}$}
\put(30,27){$\line(0,1){2}$}

\put(50,9){$\line(1,0){2}$ $\line(1,0){2}$ $\line(1,0){2}$
$\line(1,0){2}$}

\put(50,11){$\line(0,1){2}$} \put(50,15){$\line(0,1){2}$}
\put(50,19){$\line(0,1){2}$} \put(50,23){$\line(0,1){2}$}
\put(50,27){$\line(0,1){2}$} \put(50,31){$\line(0,1){2}$}

\put(27,34){$\line(1,0){2}$ $\line(1,0){2}$ $\line(1,0){2}$
$\line(1,0){2}$ $\line(1,0){2}$ $\line(1,0){2}$}


\put(70,15){$\line(0,1){2}$} \put(70,11){$\line(0,1){2}$}
\put(70,19){$\line(0,1){2}$} \put(70,23){$\line(0,1){2}$}
\put(70,27){$\line(0,1){2}$} \put(70,31){$\line(0,1){2}$}
\put(70,35){$\line(0,1){2}$}

\put(47,38){$\line(1,0){2}$ $\line(1,0){2}$ $\line(1,0){2}$
$\line(1,0){2}$ $\line(1,0){2}$ $\line(1,0){2}$}

\put(47,40){$\line(0,1){2}$} \put(47,44){$\line(0,1){2}$}
\put(47,48){$\line(0,1){2}$} \put(47,52){$\line(0,1){2}$}
\put(47,56){$\line(0,1){2}$}


\put(64,40){$\line(1,0){2}$ $\line(1,0){2}$  $\line(1,0){2}$
$\line(1,0){2}$   $\line(1,0){2}$  $\line(1,0){2}$}

\put(64,42){$\line(0,1){2}$} \put(64,46){$\line(0,1){2}$}
\put(64,50){$\line(0,1){2}$} \put(64,54){$\line(0,1){2}$}

\put(52,56){$\line(1,0){2}$ $\line(1,0){2}$  $\line(1,0){2}$}

\put(95,30){$\line(1,0){2}$ $\line(1,0){2}$  $\line(1,0){2}$}

\put(95,33){$\line(0,1){2}$} \put(95,37){$\line(0,1){2}$}
\put(95,41){$\line(0,1){2}$} \put(95,45){$\line(0,1){2}$}
\put(95,49){$\line(0,1){2}$}

\put(72,52){$\line(1,0){2}$ $\line(1,0){2}$  $\line(1,0){2}$
$\line(1,0){2}$   $\line(1,0){2}$  $\line(1,0){2}$}

\put(115,32){$\line(1,0){2}$ $\line(1,0){2}$  $\line(1,0){2}$}

\put(115,33){$\line(0,1){2}$} \put(115,37){$\line(0,1){2}$}
\put(115,41){$\line(0,1){2}$} \put(115,45){$\line(0,1){2}$}
\put(115,49){$\line(0,1){2}$} \put(115,53){$\line(0,1){2}$}

\put(95,55){$\line(1,0){2}$ $\line(1,0){2}$  $\line(1,0){2}$
$\line(1,0){2}$   $\line(1,0){2}$}

\put(95,57){$\line(0,1){2}$} \put(95,61){$\line(0,1){2}$}
\put(95,65){$\line(0,1){2}$} \put(95,69){$\line(0,1){2}$}
\put(95,74){$\line(0,1){2}$}

\put(130,33){$\line(0,1){2}$} \put(130,37){$\line(0,1){2}$}
\put(130,41){$\line(0,1){2}$} \put(130,45){$\line(0,1){2}$}
\put(130,49){$\line(0,1){2}$} \put(130,53){$\line(0,1){2}$}

\put(108,57){$\line(1,0){2}$ $\line(1,0){2}$ $\line(1,0){2}$
$\line(1,0){2}$ $\line(1,0){2}$ $\line(1,0){2}$}

\put(108,65){$\line(0,1){2}$} \put(108,69){$\line(0,1){2}$}
\put(108,74){$\line(0,1){2}$} \put(108,59){$\line(0,1){2}$}
\put(108,79){$\line(0,1){2}$}

\put(88,82){$\line(1,0){2}$ $\line(1,0){2}$  $\line(1,0){2}$
$\line(1,0){2}$   $\line(1,0){2}$}


\put(135,63){$\line(0,1){2}$} \put(135,67){$\line(0,1){2}$}
\put(135,59){$\line(0,1){2}$} \put(135,55){$\line(0,1){2}$}
\put(135,51){$\line(0,1){2}$} \put(135,71){$\line(0,1){2}$}

\put(115,73){$\line(1,0){2}$ $\line(1,0){2}$ $\line(1,0){2}$
$\line(1,0){2}$ $\line(1,0){2}$}


\put(135,82){$\line(1,0){2}$ $\line(1,0){2}$ $\line(1,0){2}$
$\line(1,0){2}$ $\line(1,0){2}$}

\put(148,73){$\bullet$}

\end{picture}
\end{center}
\caption{Snakes of the skew partition
$(7,6,6,3)/(3,1)$}\label{snake}
\end{figure}

We follow the terminology of Stanley \cite{S2} on snake sequences
and interval sets, which are helpful notions for the enumeration
of the minimal border strip decompositions of a skew partition
$\lambda/\mu$. Let us consider the bottom-right boundary lattice
path with steps $(0,1)$ or $(1,0)$ from the bottom-leftmost point
of the diagram of $\lambda/\mu$ to the top-rightmost point. We
regard this path as a sequence of edges
$e_1,\,e_2,\,\ldots,\,e_k$. For an edge $e$ in this path we define
a  subset $S_{e}$ of squares of $\lambda/\mu$, called a
\emph{snake}. If there exists no square having $e$ as an edge,
then we have the set $S_e=\emptyset$. Let $(i,j)$ be the unique
square of $\lambda/\mu$ having $e$ as an edge. If $e$ is
horizontal, then we define
\begin{equation}\label{right-snake}
S_e=\lambda/\mu\cap\{(i,j),\,(i-1,j),\,(i-1,j-1),\,(i-2,j-1),\,(i-2,j-2),\,\ldots\}.
\end{equation}
If $e$ is vertical, we then define
\begin{equation}\label{left-snake}
S_e=\lambda/\mu\cap\{(i,j),\,(i,j-1),\,(i-1,j-1),\,(i-1,j-2),\,(i-2,j-2),\,\ldots\}.
\end{equation}
For example, the nonempty snakes of the skew shape
$(7,6,6,3)/(3,1)$ are shown in Figure \ref{snake}, and the two
snakes with just one square are shown with a single bullet. The
\emph{length} $\ell(S)$ of a snake $S$ is defined to be one less
than its number of squares. For an empty snake $S$, let
$\ell(S)=-1$. A \emph{right snake} is a snake of even length and
of the form \eqref{right-snake}, and a \emph{left snake} is a
snake of even length and of the  form \eqref{left-snake}. From the
boundary lattice path we obtain a sequence of snakes:
$(S_{e_1},\,S_{e_2},\,\ldots,\,S_{e_k})$. The \emph{snake
sequence} of $\lambda/\mu$, denoted ${\rm SS}(\lambda/\mu)$, is
defined by replacing a left snake of length $2m$ with the symbol
$L_m$ in the sequence $(S_{e_1},\,S_{e_2},\,\ldots,\,S_{e_k})$,
replacing a right snake of length $2m$ with $R_m$, and replacing a
snake of odd length with $O$. From Figure \ref{snake}, we see that
$${\rm SS}((7,6,6,3)/(3,1))=L_0 L_1 O O O O L_2 R_2 R_1 O R_0.$$

Let ${\rm rank}(\lambda/\mu)=r$, and let ${\rm
SS}(\lambda/\mu)=q_1q_2\cdots q_k$. An \emph{interval set}
$\mathcal{I}$ of $\lambda/\mu$ is defined to be a collection of
$r$ ordered pairs $\{(u_1,v_1),\,(u_2,v_2),\,\ldots,\,(u_r,v_r)\}$
such that
\begin{enumerate}
\item $u_i\neq u_j$ and $v_i\neq v_j$ for $1\leq i<j\leq r$.

\item $1\leq u_i<v_i\leq k$ and $u_i\neq v_j$ for $1\leq i,j\leq
r$.

\item $q_{u_i}=L_s$ and $q_{v_i}=R_{s'}$ for some $s$ and $s'$
(depending on $i$).
\end{enumerate}

Let ${\rm cr}(\mathcal{I})$ denote the number of crossings of
$\mathcal{I}$, i.e., the number of pairs $(i,j)$ for which
$u_i<u_j<v_i<v_j$. According to \cite[Proposition 4.3]{S2}, there
exists a unique interval set
$\mathcal{I}_0=\{(w_1,y_1),\,(w_2,y_2),\,\ldots,\,(w_r,y_r)\}$
such that ${\rm cr}(\mathcal{I}_0)=0$. From \cite{S2},  we see
that ${\rm SS}(\lambda/\mu)$ has exactly $r$ left snakes and $r$
right snakes. For an interval set
$\mathcal{I}=\{(u_1,v_1),\,(u_2,v_2),\,\ldots,\,(u_r,v_r)\}$,
 we
may impose a linear order $u_1<u_2<\cdots<u_r$ on its elements.
Then there exists a unique permutation $\sigma$ relative to
$\mathcal{I}_0$ such that for each $i$
\begin{equation}
u_i=w_i \mbox{ and } v_i=y_{\sigma_i}.
\end{equation}
Thus, each interval set $\mathcal{I}$ is associated to a
permutation $\sigma_I$, which we call the \emph{interval
permutation} of $\mathcal{I}$ with respect to $\mathcal{I}_0$.
Given a permutation $\sigma$, let $\rm{inv}(\sigma)$ denote the
number of inversions of $\sigma$, i.e., the number of pairs
$(i,j)$ satisfying $i<j$ but $\sigma_i>\sigma_j$.

\begin{prop}\label{key-le}
Given a skew partition $\lambda/\mu$ and an interval set
$\mathcal{I}$ of $\lambda/\mu$, let $\sigma_{\mathcal{I}}$ be the
interval permutation with respect to $\mathcal{I}_0$. Then we have
\begin{equation}
{\rm cr}(\mathcal{I})\equiv \rm{inv}(\sigma_{\mathcal{I}})\
(\rm{mod}\ 2).
\end{equation}
\end{prop}

\proof First we give a geometric representation of ${\rm
cr}(\mathcal{I})$. For each interval $(u_i,\,v_i)$ of
$\mathcal{I}$ we draw an arc on top of ${\rm SS}(\lambda/\mu)$
which connects two snakes $q_{u_i}$ and $q_{v_i}$. For a given
pair $(i,j)$ with $i<j$, the two arcs $(u_i,\,v_i)$ and
$(u_j,\,v_j)$ are said to be noncrossing if $u_i<u_j<v_j<v_i$. In
this terminology ${\rm cr}(\mathcal{I})$ equals the number of
crossings.

To determine the inversions of $\sigma_{\mathcal{I}}$, we replace
$q_{w_i}$ by $F_i$ and $q_{y_i}$ by $G_i$ in ${\rm
SS}(\lambda/\mu)$ for each $i$. Clearly, $\sigma_{\mathcal{I}}$ is
a bijection from $\{F_1,\,F_2,\,\ldots,\,F_r\}$ to
$\{G_1,\,G_2,\,\ldots,\,G_r\}$. We now represent the snakes of
${\rm SS}(\lambda/\mu)$ with respect to the order $F_1,F_2,\cdots,
F_r,G_r,G_{r-1},\cdots, G_1$ by moving $G_1$ to the right of the
rightmost element if $G_1$ itself is not the rightmost element,
and repeating this process until we achieve the desired order.  It
follows that $\rm{inv}(\sigma_{\mathcal{I}})$ equals the number of
crossings in the above representation. Note that at each step of
moving $G_i$ to the proper position, the number of crossings in
the diagram can only change by an even number. This completes the
proof.
  \qed

\begin{figure}[h,t]
\begin{center}
\setlength{\unitlength}{2pt}
\begin{picture}(160,35)
\put(40,5){{\rm
$L_0\,O\,L_1\,L_2\,R_2\,O\,O\,L_2\,R_2\,O\,R_1\,R_0$}}
\qbezier(43,13)(80,35)(117,13) \qbezier(55,13)(82,30)(109,13)
\qbezier(62,13)(66,17)(70,13) \qbezier(89,13)(93,17)(97, 13)
\end{picture}
\end{center}
\caption{Parenthesization of the snake sequence ${\rm
SS}((8,8,7,4)/(4,1,1))$}\label{figure-snake}
\end{figure}
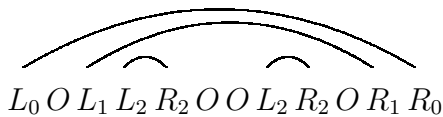

For example, let $\lambda/\mu=(8,8,7,4)/(4,1,1)$. Figure
\ref{figure-snake} shows the snake sequence ${\rm
SS}((8,8,7,4)/(4,1,1))$, from which we see that
\begin{equation*}
\mathcal{I}_0=\{(1, 12), (3, 11), (4, 5), (8, 9)\}.
\end{equation*}

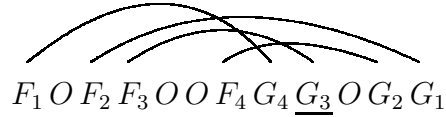
\begin{figure}[h,t]
\begin{center}
\setlength{\unitlength}{2pt}
\begin{picture}(160,35)
\put(40,5){{\rm
$F_1\,O\,F_2\,F_3\,O\,O\,F_4\,G_4\,\underline{G_3}\,O\,G_2\,G_1$}}
\qbezier(43,13)(70,35)(89,13) \qbezier(55,13)(81,30)(117,13)
\qbezier(62,13)(80,25)(97, 13) \qbezier(80,13)(90,20)(109,13)
\end{picture}
\end{center}
\caption{After moving $G_3$}\label{figure-a}
\end{figure}

Let us illustrate the proof of Proposition \ref{key-le} by the
example  \break $\mathcal{I}=\{(1, 9), (3, 12), (4, 5), (8,
11)\}$, for which we have $\sigma_{\mathcal{I}}=[4, 1, 3, 2 ]$.
The crossings of $\mathcal{I}$ are shown in Figure \ref{figure-b},
where we relabel the snakes as described in the proof.  Figure
\ref{figure-a} demonstrates the diagram after moving $G_3$ which
has two more crossings. It is evident that
\begin{equation}
{\rm cr}(\mathcal{I})=2, \quad {\rm inv}(\pi)=4, \quad {\rm
cr}(\mathcal{I})\equiv \rm{inv}(\pi)\ (\rm{mod}\ 2).
\end{equation}
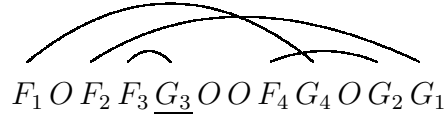
\begin{figure}[h,t]
\begin{center}
\setlength{\unitlength}{2pt}
\begin{picture}(160,35)
\put(40,5){{\rm
$F_1\,O\,F_2\,F_3\,\underline{G_3}\,O\,O\,F_4\,G_4\,O\,G_2\,G_1$}}

\qbezier(43,13)(70,35)(97, 13) \qbezier(55,13)(81,30)(117,13)
\qbezier(62,13)(66,17)(70,13) \qbezier(89,13)(99,17)(109,13)

\end{picture}
\end{center}
\caption{Before moving $G_3$}\label{figure-b}
\end{figure}

\section{Minimal border strip decompositions}

\begin{figure}[h,t]
\setlength{\unitlength}{2pt}
\begin{center}
\begin{picture}(186,100)
\put(80,85){\line(1,0){60}} \put(60,65){\line(1,0){80}}
\put(60,45){\line(1,0){60}} \put(40,25){\line(1,0){60}}
\put(40,5){\line(1,0){40}}

\put(40,5){\line(0,1){20}} \put(60,5){\line(0,1){60}}
\put(80,5){\line(0,1){80}} \put(100,25){\line(0,1){60}}
\put(120,45){\line(0,1){40}} \put(140,65){\line(0,1){20}}

\put(142, 73){{1}}\put(130, 55){0}\put(122, 53){{1}} \put(110,
35){0}\put(102, 33){{1}} \put(90, 15){0}\put(82, 13){{1}} \put(70,
-5){0} \put(50, -5){0}

\put(75, 73){{1}}\put(130, 90){0}\put(55, 53){{1}} \put(110,
90){0}\put(55, 33){{1}} \put(90, 90){0}\put(35, 13){{1}} \put(70,
67){0} \put(50, 27){0}
\end{picture}
\end{center}
\caption{Constructing the reduced code of
$(5,4,3,2)/(2,1,1)$}\label{bound}
\end{figure}
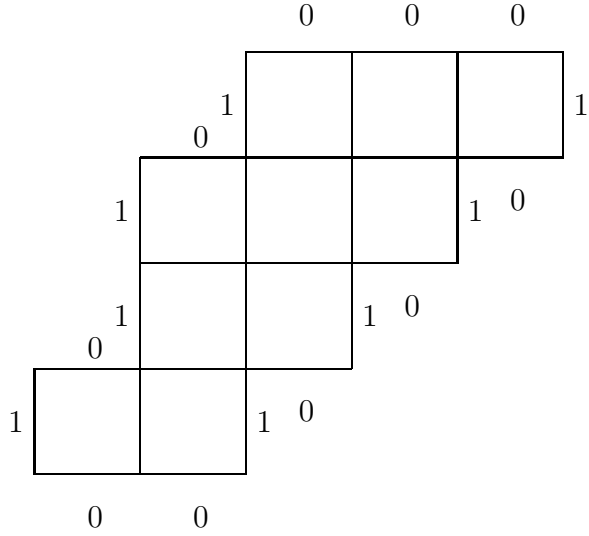

We recall the notion of the \emph{reduced code} of a skew
partition $\lambda/\mu$, denoted ${\rm c}(\lambda/\mu)$. The
reduced code ${\rm c}(\lambda/\mu)$ is also known as the
\emph{partition sequence} of $\lambda/\mu$ \cite{B1, B2}. Consider
the two boundary lattice paths of the diagram of $\lambda/\mu$
with steps $(0,1)$ or $(1,0)$ from the bottom-leftmost point to
the top-rightmost point. Replacing each step $(0,1)$ by $1$ and
each step $(1,0)$ by $0$,  we obtain two binary sequences by
reading the lattice paths from the bottom-left corner to the
top-right corner. Denote the top-left binary sequence by
$f_1,\,f_2,\,\ldots,\,f_k$, and the bottom-right binary sequence
by $g_1,\,g_2,\,\ldots,\,g_k$. The reduced code ${\rm
c}(\lambda/\mu)$ is defined by the two-line array
$$
\begin{array}{cccc}
f_1 & f_2 & \cdots & f_k\\
g_1 & g_2 & \cdots & g_k
\end{array}.
$$
The reduced code of the skew partition  $(5,4,3,2)/(2,1,1)$ in
Figure \ref{bound} is
$$
\begin{array}{ccccccccc}1 & 0 & 1 & 1 & 0 &
1 & 0 & 0 & 0\\
0 & 0 & 1 & 0 & 1 & 0 & 1 & 0 & 1
\end{array}.
$$

A \emph{diagonal} with content $j$ of $\lambda/\mu$, denoted
$d_j(\lambda/\mu)$, is the set of all the squares in $\lambda/\mu$
having content $j$. Suppose that the length of ${\rm
c}(\lambda/\mu)$ is $k$. It is obvious that $\lambda/\mu$ has
$k-1$ diagonals. Let $\epsilon$ be the smallest content of
$\lambda/\mu$. For each $i:1\leq i\leq k-1$, we put the diagonal
$d_{\epsilon+i-1}$ between the $i$-th column and $(i+1)$-th column
of ${\rm c}(\lambda/\mu)$. Then we obtain a connection between the
diagonals of $\lambda/\mu$ and the reduced code ${\rm
c}(\lambda/\mu)$.

Recall that a skew partition $\lambda/\mu$ is said to be
\emph{connected} if the interior of the Young diagram of
$\lambda/\mu$ is a connected set. A \emph{border strip} is a
connected skew partition with no $2\times 2$ square. Define the
size of a border strip $B$ as the number of squares of $B$, and
define the \emph{height} $ht(B)$ of $B$ as one less than its
number of rows. We say that $B\subset \lambda/\mu$ is a border
strip of $\lambda/\mu$ if $\lambda/\mu-B$ is a skew partition
$\nu/\mu$. A border strip $B$ of $\lambda/\mu$ is said to be
\emph{maximal} if there does not exist another border strip
$B'\subset \lambda/\mu$ such that $B\subset B'$. A \emph{border
strip decomposition} \cite{S1} of $\lambda/\mu$ is a partition of
the squares of $\lambda/\mu$ into pairwise disjoint border strips.
A \emph{greedy border strip decomposition} of $\lambda/\mu$ is
obtained by successively removing the maximal border strip from
$\lambda/\mu$. A border strip decomposition is \emph{minimal} if
there does not exist a border strip decomposition with a fewer
number of border strips.

 Stanley \cite[Proposition 2.2]{S2} has shown that the
rank of a skew partition $\lambda/\mu$ is equal to the number of
border strips in a minimal border strip decomposition of
$\lambda/\mu$, and it is also equal to the number of ${1 \atop 0}$
columns of ${\rm c}(\lambda/\mu)$. As a consequence, a greedy
border strip decomposition is minimal, because when we
successively remove the maximal border strips from $\lambda/\mu$ a
column ${1\atop 0}$ of ${\rm c}(\lambda/\mu)$ changes into
${1\atop 1}$ and a column ${0\atop 1}$ changes into ${0\atop 0}$.

Suppose that ${\rm rank}(\lambda/\mu)=r$. Given a minimal border
strip decomposition $\mathbf{D}=\{B_1,\,B_2,\,\ldots,\,B_r\}$ of
$\lambda/\mu$, let
$$P_{\mathbf{D}}=\{\tau({\rm init}(B_1)),\,\tau({\rm init}(B_2)),\,\ldots,
\tau({\rm init}(B_r))\}$$ and
$$Q_{\mathbf{D}}=\{\tau({\rm fin}(B_1)),\,\tau({\rm fin}(B_2)),\,\ldots,
\tau({\rm fin}(B_r))\},$$ where ${\rm init}(B_i)$ is the lower
left-hand square of $B_i$ and ${\rm fin}(B_i)$ is the upper
right-hand square. The following proposition shows that
$P_{\mathbf{D}}$ and $Q_{\mathbf{D}}$ are independent of the
minimal border strip decomposition $\mathbf{D}$.

\begin{prop} \label{cont-bord} Let
$\mathcal{I}_0=\{(w_1,y_1),\,(w_2,y_2),\,\ldots,\,(w_r,y_r)\}$ be
the interval set of $\lambda/\mu$ with ${\rm
cr}(\mathcal{I}_0)=0$. Let $\epsilon$ be the smallest value among
the contents of the squares of $\lambda/\mu$. Let $\mathbf{D}$ be
a minimal border strip decomposition of $\lambda/\mu$. Then we
have
\begin{equation}
P_{\mathbf{D}}=\{\epsilon+w_i-1\, |\, 1\leq i\leq r\} \mbox{ and }
Q_{\mathbf{D}}=\{\epsilon+y_i-2\, |\, 1\leq i\leq r\}.
\end{equation}
\end{prop}

\proof By \cite[Proposition 2.1]{S2}, we see that the operation of
removing a border strip $B$ of size $p$ from $\lambda/\mu$
corresponds to the operation of choosing $i$ with the $i$-th
column being ${1\atop 0}$ and the $(i+p)$-th column being ${0\atop
1}$, and then replacing the $i$-th column with ${1\atop 1}$ and
the $(i+p)$-th column with ${0\atop 0}$. Moreover, the lower
left-hand square of $B$ lies on the diagonal $d_{i}$, and the
upper right-hand square of $B$ lies on the diagonal $d_{i+p-1}$.
Therefore
\[\tau({\rm init}(B))=\epsilon+i-1 \mbox{ and }
\tau({\rm fin}(B))=\epsilon+i+p-2.\]
 It follows that $P_{\mathbf{D}}$
and $Q_{\mathbf{D}}$ are determined by the indices of the columns
${1\atop 0}$ and ${0\atop 1}$ of ${\rm c}(\lambda/\mu)$
respectively. Since $\{w_i\}$ is  the set of indices of columns
${1\atop 0}$ of ${\rm c}(\lambda/\mu)$, and $\{y_i\}$ is the set
of indices of ${0\atop 1}$, we get  the desired assertion. \qed

\section{Giambelli-type determinantal formulas }

In this section, we obtain a determinantal formula for the
quantity given by Stanley based on the Giambelli-type formula for
skew Schur functions. Let $\lambda/\mu$ be a skew diagram. A
border strip decomposition of $\lambda/\mu$ is said to be an
\emph{outside decomposition} if every strip in the decomposition
has an initial square on the left or bottom perimeter of the
diagram and a terminal square on the right or top perimeter, see
Figure \ref{bd1}. It is obvious that a greedy border strip
decomposition of $\lambda/\mu$ is an outside decomposition.

\begin{figure}[h,t]
\begin{center}
\setlength{\unitlength}{20pt}
\begin{picture}(16,4)

\put(1,1){\line(1,0){2.1}} \put(1,1.7){\line(1,0){3.5}}
\put(1.7,2.4){\line(1,0){3.5}} \put(1.7,3.1){\line(1,0){3.5}}
\put(2.4,3.8){\line(1,0){2.8}}

\put(1,1){\line(0,1){0.7}}\put(1.7,1){\line(0,1){2.1}}
\put(2.4,1){\line(0,1){2.8}}\put(3.1,1){\line(0,1){2.8}}
\put(3.8,1.7){\line(0,1){2.1}}\put(4.5,1.7){\line(0,1){2.1}}
\put(5.2,2.4){\line(0,1){1.4}}

\multiput(1.35,1.35)(0.35,0){4}{\line(1,0){0.175}}
\multiput(2.75,3.45)(0.35,0){6}{\line(1,0){0.175}}
\multiput(4.15,2.75)(0.35,0){2}{\line(1,0){0.175}}

\multiput(2.75,1.35)(0,0.35){6}{\line(0,1){0.175}}
\multiput(2.05,1.875)(0,0.35){3}{\line(0,1){0.175}}
\multiput(3.45,1.875)(0,0.35){3}{\line(0,1){0.175}}
\multiput(4.15,1.875)(0,0.35){3}{\line(0,1){0.175}}

\put(10,1){\line(1,0){2.1}} \put(10,1.7){\line(1,0){3.5}}
\put(10.7,2.4){\line(1,0){3.5}} \put(10.7,3.1){\line(1,0){3.5}}
\put(11.4,3.8){\line(1,0){2.8}}

\put(10,1){\line(0,1){0.7}}\put(10.7,1){\line(0,1){2.1}}
\put(11.4,1){\line(0,1){2.8}}\put(12.1,1){\line(0,1){2.8}}
\put(12.8,1.7){\line(0,1){2.1}}\put(13.5,1.7){\line(0,1){2.1}}
\put(14.2,2.4){\line(0,1){1.4}}

\multiput(10.35,1.35)(0.35,0){2}{\line(1,0){0.175}}
\multiput(11.05,2.75)(0.35,0){2}{\line(1,0){0.175}}
\multiput(12.45,3.45)(0.35,0){4}{\line(1,0){0.175}}
\multiput(11.75,2.05)(0.35,0){2}{\line(1,0){0.175}}
\multiput(13.15,2.75)(0.35,0){2}{\line(1,0){0.175}}

\multiput(11.05,1.35)(0,0.35){4}{\line(0,1){0.175}}
\multiput(11.75,2.75)(0,0.35){3}{\line(0,1){0.175}}
\multiput(12.45,2.05)(0,0.35){4}{\line(0,1){0.175}}
\multiput(11.75,1.35)(0,0.35){2}{\line(0,1){0.175}}
\multiput(13.15,2.05)(0,0.35){2}{\line(0,1){0.175}}

\put(-0.5,0){a. A border strip decomposition} \put(9,0){b. An
outside decomposition}

\end{picture}
\end{center}
\caption{Border strip decompositions} \label{bd1}
\end{figure}

The notion of the cutting strip of an outside decomposition is
introduced by Chen, Yan and Yang \cite{CYY}, which is used  to
give a transformation theorem on the Giambelli-type determinantal
formulas for the skew Schur function.

We proceed to construct a cutting strip for an edgewise connected
skew partition $\lambda/\mu$. Suppose that $\lambda/\mu$ has $k$
diagonals. The cutting strip of an outside decomposition is
defined to be a border strip of length $k$. Given an outside
decomposition, we may assign a direction to each square in the
diagram. Starting with the bottom-left corner of a strip, we say
that a square of a strip has up direction (resp. right direction)
if the next square in the strip lies on its top (resp. to its
right). Notice that the strips in any outside decomposition of
$\lambda/\mu$ are nested in the sense that the squares in the same
diagonal of $\lambda/\mu$ all have up direction or all have right
direction. Based on this property, the cutting strip $\phi$ of an
outside decomposition $\mathbf{D}$ of $\lambda/\mu$ is defined as
follows: for $i=1,\,2,\,\ldots,\,k-1$ the $i$-th square in $\phi$
keeps the same direction as the $i$-th diagonal of $\lambda/\mu$
with respect to $\mathbf{D}$. For any two integers $p,q$ a strip
$[p,q]$ is defined by the following rule: if $p\leq q$, then let
$[p,q]$ be the segment of $\phi$ from the square with content $p$
to the square with content $q$; if $p=q+1$, then let $[p,q]$ be
the empty strip; if $p>q+1$, then $[p,q]$ is undefined. Using the
above notation, Hamel and Goulden's theorem on the Giambelli-type
formulas for the skew Schur function can be formulated as follows.

\begin{theo}[{\cite[Theorem 3.1]{HG}}]\label{schur-dec}
For an outside decomposition $\mathbf{D}$ with $k$ border strips
$B_1,\,B_2,\,\ldots,\,B_k$, we have
\begin{equation}
s_{\lambda/\mu}=\det\left(s_{[\tau({\rm init}(B_i)),\tau({\rm
fin}(B_j))]}\right)_{i,j=1}^{k}.
\end{equation}
\end{theo}
By choosing the outside decomposition whose border strips are the
rows of the diagram of $\lambda/\mu$ in the above theorem, we
obtain the Jacobi-Trudi identity for the skew Schur function,
which states that
\begin{equation}
s_{\lambda/\mu}=\det\left(h_{\lambda_i-\mu_j-i+j}\right)_{i,j=1}^{\ell(\lambda)},
\end{equation}
where $h_k$ denotes the $k$-th complete symmetric function,
$h_0=1$ and $h_k=0$ for $k<0$.

Let $y(\lambda/\mu)=(t^{-{\rm
rank}(\lambda/\mu)}s_{\lambda/\mu}(1^t))_{t=0}$. The zrank
conjecture says that $y(\lambda/\mu)\neq 0$ for any skew partition
$\lambda/\mu$. Now we give the evaluation of $y(\lambda/\mu)$ by
using Theorem \ref{schur-dec}. First we consider the case when
$\lambda/\mu$ is a border strip. In this case we have ${\rm
rank}(\lambda/\mu)=1$, $\mu_i=\lambda_{i+1}-1$ for $i\leq
\ell(\lambda)-1$ and $\mu_{\ell(\lambda)}=0$. From the
Jacobi-Trudi identity one easily deduces the following lemma.

\begin{lemm}\label{ribbon-lemm} For a border strip $\lambda/\mu$
we have
\begin{equation}\label{ribbon-eq}
y(\lambda/\mu)=\frac{(-1)^{\ell(\lambda)+1}}{\lambda_1+\ell(\lambda)-1}.
\end{equation}
\end{lemm}

In order to compute $y(\lambda/\mu)$ for a general skew partition
$\lambda/\mu$, we need to consider the greedy border strip
decomposition
 $\mathbf{D_0}$ of $\lambda/\mu$. Suppose that ${\rm rank}(\lambda/\mu)=r$.
It follows that $\mathbf{D_0}$ has $r$ border strips. We may apply
Theorem \ref{schur-dec} to $\mathbf{D_0}$ because it is also an
outside decomposition. Furthermore, we may impose a canonical
order on the strips $B_1,\,B_2,\,\ldots,\,B_r$ of $\mathbf{D_0}$
by the contents of their lower left-hand squares such that
$\tau({\rm init}(B_i))< \tau({\rm init}(B_{i+1}))$ for $i<r$.
Since the sum  of the heights of border strips in $\mathbf{D}_0$
is uniquely determined by the shape $\lambda/\mu$,  one sees that
$$z(\lambda/\mu)=ht(B_1)+ht(B_2)+\cdots+ht(B_r)$$
is well defined. Let
$\mathcal{I}_0=\{(w_1,y_1),\,(w_2,y_2),\,\ldots,\,(w_r,y_r)\}$ be
the interval set of $\lambda/\mu$ with ${\rm
cr}(\mathcal{I}_0)=0$. By Proposition \ref{cont-bord} and the
properties of $\mathbf{D}_0$ and $\mathcal{I}_0$, we obtain that
\begin{equation}\label{con-int}
\tau({\rm init}(B_i))=\epsilon+w_i-1 \mbox{ and } \tau({\rm
fin}(B_i))=\epsilon+y_i-2,
\end{equation}
where $\epsilon$ is the smallest value among the contents of the
squares of $\lambda/\mu$.

The following theorem gives a determinantal formula for
$y(\lambda/\mu)$ based on a matrix related to the Cauchy matrix.

\begin{theo}\label{th-main}Let $\lambda/\mu$ be a skew partition
with ${\rm rank}(\lambda/\mu)=r$, and let $\mathcal{I}_0$ be the
noncrossing interval set
$\{(w_1,y_1),\,(w_2,y_2),\,\ldots,\,(w_r,y_r)\}$ of $\lambda/\mu$.
Then we have
\begin{equation}\label{th-main-eq}
y(\lambda/\mu)=(-1)^{z(\lambda/\mu)}\det(d_{ij})_{i,j=1}^r,
\end{equation}
where
$$
d_{ij}=\left\{
\begin{array}{ll}
\displaystyle\frac{1}{y_j-w_i}, & \mbox{if } y_j>w_i\\[12pt]
0, & \mbox{if } y_j<w_i
\end{array}
\right.
$$
\end{theo}

\proof Take the greedy outside decomposition
$\textbf{D}_0=\{B_1,\,B_2,\,\ldots,\,B_r\}$ of $\lambda/\mu$, and
let $\phi_0$ be the cutting strip corresponding to $\textbf{D}_0$.
By Theorem \ref{schur-dec} we have
\begin{equation}
s_{\lambda/\mu}=\det\left(s_{[\tau({\rm init}(B_i)),\tau({\rm
fin}(B_j))]}\right)_{i,j=1}^{r}.
\end{equation}
Suppose that the square with content $\tau({\rm init}(B_i))$ lies
in the $p_i$-th row of $\phi_0$, and the square with content
$\tau({\rm fin}(B_j))$ lies in the $q_j$-th row. Applying Lemma
\ref{ribbon-lemm}, we get
\begin{equation}\label{j1}
(t^{-1}s_{[\tau({\rm init}(B_i)),\tau({\rm
fin}(B_j))]})_{t=0}=\displaystyle\frac{(-1)^{{p_i}-{q_j}}}{\tau({\rm
fin}(B_j))+1-\tau({\rm init}(B_i))}
\end{equation}
if $[\tau({\rm init}(B_i)),\tau({\rm fin}(B_j))]$ is a substrip of
$\phi_0$. Otherwise, the above entry is set $0$. Note that
$[\tau({\rm init}(B_i)),\tau({\rm fin}(B_j))]$ cannot be an empty
strip for the greedy border strip decomposition. Using
\eqref{con-int} we may write \eqref{j1} as
\begin{equation}
(t^{-1}s_{[\tau({\rm init}(B_i)),\tau({\rm
fin}(B_j))]})_{t=0}=\displaystyle\frac{(-1)^{{p_i}-{q_j}}}{y_j-w_i}
\end{equation}
for $y_j>w_i$, or $0$ for $y_j<w_i$. Thus, we have
\begin{equation}
y(\lambda/\mu)=(t^{-r}s_{\lambda/\mu}(1^t))_{t=0}=\det\left((t^{-1}s_{[\tau({\rm
init}(B_i)),\tau({\rm fin}(B_j))]})_{t=0}\right)_{i,j=1}^{r}.
\end{equation}
Extracting the signs from the determinant, we obtain
$$y(\lambda/\mu)=(-1)^{({p_1}+\cdots+{p_r})-({q_1}+\cdots+{q_r})}\det(d_{ij})_{i,j=1}^r=(-1)^{z(\lambda/\mu)}\det(d_{ij})_{i,j=1}^r.$$
This completes the proof. \qed

\noindent \textbf{Remark.} Stanley \cite{stanleyprv} pointed out
that one can also get a matrix for $y(\lambda/\mu)$ by taking the
Jacobi-Trudi matrix (the matrix appearing in the Jacobi-Trudi
determinant formula of $s_{\lambda/\mu}$) for the skew Schur
function $s_{\lambda/\mu}$, and deleting all rows and columns that
contain a $1$, and then substituting $1/i$ for $h_i$. This matrix
coincides with the matrix $(d_{ij})_{i,j=1}^r$ defined in
\eqref{th-main-eq}, subject to permutations of rows and columns.
This fact can be verified by using the transformation formula in
\cite{CYY}.

From Theorem \ref{th-main} and Proposition \ref{key-le} one can
recover the following expansion formula of Stanley \cite[Equation
(30)]{S2}.

\begin{coro}\label{coro-main} We have
\begin{equation}\label{eq-comb-v}
y(\lambda/\mu)=(-1)^{z(\lambda/\mu)}\sum_{\mathcal{I}=\{(u_1, v_1
), \ldots, (u_r, v_r)\}}\frac{(-1)^{{\rm
cr}(\mathcal{I})}}{\prod_{i=1}^r(v_i-u_i)},
\end{equation}
summed over all interval sets $\mathcal{I}$ of $\lambda/\mu$.
\end{coro}

\section{An equivalent description of the zrank conjecture}

We begin this section with the definition of a restricted Cauchy
matrix. Let $a=(a_1, \ldots, a_n)$ and $b=(b_1, \ldots, b_n)$ be
two integer sequences. Suppose that $a$ is strictly decreasing and
$b$ is strictly increasing, and for any $i,j$ we have
$a_i>b_{n+1-i}$ and $a_i\neq b_j$. We define a matrix
$C(a,b)=(c_{ij})_{i,j=1}^n$ by setting
$$
c_{ij}=\left\{
\begin{array}{ll}
{\displaystyle \frac{1}{a_i-b_j}}, & \mbox{ if $a_i>b_j$}\\[12pt]
0, & \mbox{ if $a_i<b_j$}
\end{array}
\right..
$$

\begin{deff}
A matrix $M$ is called a \emph{restricted Cauchy matrix} if there
exist two integer sequences $a$ and $b$ satisfying the above
conditions such that $M=C(a,b)$.
\end{deff}

For a matrix $M$ we say it is \emph{singular} if $\det(M)=0$; or
\emph{nonsingular}, otherwise. We now come to the main result of
this paper.

\begin{theo}\label{eq-des} The following two statements are equivalent:

(i) The zrank conjecture is true for any  skew partition.

(ii) Any restricted Cauchy matrix is nonsingular.
\end{theo}

\proof Suppose  that (ii) is true. For a skew partition
$\lambda/\mu$, consider the noncrossing interval set
$\mathcal{I}_0=\{(w_1,y_1),\,(w_2,y_2),\,\ldots,\,(w_r,y_r)\}$ of
$\lambda/\mu$. Clearly, $w_i\neq y_j$ for $1\leq i,j\leq r$. Let
$w=(w_1',\,w_2',\ldots,\,w_r')$ be  the rearrangement of
$(w_1,\,w_2,\ldots,\,w_r)$ in increasing order, and let
$y=(y_1',\,y_2',\ldots,\,y_r')$ be the rearrangement of
$(y_1,\,y_2,\ldots,\,y_r)$ in decreasing order. For $1\leq i\leq
r$, we have $y_i'>w_{r+1-i}'$ since the number of ${1\atop 0}$
columns in the first $\ell$ columns of the reduced code ${\rm
c}(\lambda/\mu)$ is bigger than or equals to the number of
${0\atop 1}$ columns for $1\leq\ell\leq k$, where $k$ is the
length of ${\rm c}(\lambda/\mu)$. Notice that the determinant
$\det(d_{ij})_{i,j=1}^r$ appearing in \eqref{th-main-eq} is equal
to the determinant of the restricted Cauchy matrix $C(y,w)$ up to
a sign. By Theorem \ref{th-main} we see that
$$
y(\lambda/\mu)\neq 0 \Leftrightarrow \det(d_{ij})_{i,j=1}^r\neq 0
\Leftrightarrow \det(C(y,w))\neq 0.
$$
Since the matrix $C(y,w)$ is nonsingular,  we have ${\rm
rank}(\lambda/\mu)={\rm zrank}(\lambda/\mu)$.

Now we proceed to prove (ii) by assuming that (i) is true. Given a
restricted Cauchy matrix $C(a,b)$ of order $r$, without loss of
generality, we may assume that  $a$ and $b$ are sequences of
positive integers. Let $\lambda$ be the partition with
$\lambda_i=a_i-r+i$, and let $\mu$ be the partition with
$\mu_i=b_{r+1-i}-r+i$. From $a_i>b_{r+1-i}$ we may deduce
$\lambda_i>\mu_i$ for all $i$. Thus we can construct a skew
diagram $\lambda/\mu$. Observe that the Jacobi-Trudi matrix
$(h_{\lambda_i-\mu_j-i+j})$  of $s_{\lambda/\mu}$ does not have a
column containing $1$ since
$$\lambda_i-\mu_j-i+j=a_i-b_{r+1-j}\neq 0, \mbox{for $1\leq i,j\leq r$}.$$
It follows that ${\rm rank}(\lambda/\mu)=r$ from
\cite[Proposition]{S2}. Therefore, we have
$$y(\lambda/\mu)=(t^{-r}s_{\lambda/\mu}(1^t))_{t=0}=\det\left((t^{-1}h_{\lambda_i-\mu_j-i+j}(1^t))_{t=0}\right)_{i,j=1}^{r},$$
which is the determinant $\det(C(a,b))$ up to a sign. If the zrank
conjecture is true for $\lambda/\mu$, then we have
$y(\lambda/\mu)\neq 0$, implying that  $C(a,b)$ is nonsingular.
This completes the proof.  \qed

We remark that we may restrict our attention to irreducible
restricted Cauchy matrices for the verification  of the zrank
conjecture. In other words, if every irreducible restricted Cauchy
matrix is nonsigular, then every restricted Cauchy matrix is
nonsingular.

\section{Special Cases}

In this section we consider several classes of  restricted Cauchy
matrices $C(a,b)=(c_{ij})_{i,j=1}^r$ for which we can prove that
they are nonsingular.

\noindent\textbf{Class I.} For all $i,j$ we have $r_{ij}\neq 0$.

In this case, $(c_{ij})_{i,j=1}^r$ is a Cauchy matrix.  Cauchy
\cite{Muir} showed that
\begin{equation}
\det\left(\frac{1}{a_i-b_j}\right)_{i,j=1}^r=\prod_{i<j}(a_i-a_j)\prod_{i<j}(b_j-b_i)\prod_{i,j}\frac{1}{a_i-b_j}.
\end{equation}
It follows that  $$\det(c_{ij})_{i,j=1}^r>0.$$

From the proof of  {\cite[Theorem 3.2 (b)]{S2}}, we get

\begin{prop}\label{sp} For a connected skew diagram
$\lambda/\mu$, if every row of the Jacobi-Trudi matrix that
contains a $0$ also contains a $1$, then the matrix
$(d_{ij})_{i,j=1}^r$ appearing in \eqref{th-main-eq} must satisfy
that $d_{ij}\neq 0$ for all $i,j$.
\end{prop}

Theorem \ref{th-main} and Proposition \ref{sp} yield another proof
of {\cite[Theorem 3.2]{S2}} of Stanley. Some skew partitions do
not have the property stated in the above proposition, but the
matrices $(d_{ij})_{i,j=1}^r$ are Cauchy matrices. For instance,
taking $\lambda/\mu=(8,8,7,7,7,6,1)/(5,5,3,3,2)$, its Jacobi-Trudi
matrix is
$$
\displaystyle
s_{(8,8,7,7,7,6,1)/(5,5,3,3,2)}=\displaystyle\begin{vmatrix}
h_3 & h_4 & h_7 & h_8 & h_{10} & h_{13} & h_{14}\\[8pt]
h_2 & h_3 & h_6 & h_7 & h_9 & h_{12} & h_{13}\\[8pt]
1  & h_1 & h_4 & h_5 & h_7 & h_{10} & h_{11}\\[8pt]
0  & 1  & h_3 & h_4 & h_6 & h_9 & h_{10}\\[8pt]
0  & 0  & h_2 & h_3 & h_5 & h_8 & h_9\\[8pt]
0  & 0  & 1  & h_1 & h_3 & h_6 & h_7\\[8pt]
0  & 0  & 0  & 0  & 0  & 1  & h_1
\end{vmatrix}.
$$

\noindent\textbf{Class II.} For all $(i, j)\neq (r, r)$, we have
$c_{ij}\neq 0$ and $c_{rr}=0$.

Let
$$M=\prod_{i=1}^{r-1}\frac{(a_r-a_i)(b_i-b_r)}{(a_r-b_i)(a_i-b_r)}.$$
Since $b_r>a_r$, it is easy to show that $M>1$. We see that the
restricted Cauchy matrix in this case is of the following form:
$$(c_{ij})_{i,j=1}^r=
\begin{pmatrix}
\displaystyle\frac{1}{a_1-b_1} & \ldots & \displaystyle\frac{1}{a_1-b_{r-1}} & \displaystyle\frac{1}{a_1-b_r}\\
\displaystyle\frac{1}{a_2-b_1} & \ldots & \displaystyle\frac{1}{a_2-b_{r-1}} & \displaystyle\frac{1}{a_2-b_r}\\
\vdots            & \ldots & \vdots                & \vdots\\
\displaystyle\frac{1}{a_{r-1}-b_1} & \ldots & \displaystyle\frac{1}{a_{r-1}-b_{r-1}} & \displaystyle\frac{1}{a_{r-1}-b_r}\\
\displaystyle\frac{1}{a_r-b_1} & \ldots &
\displaystyle\frac{1}{a_r-b_{r-1}} & 0
\end{pmatrix}.$$
Then we have
\begin{eqnarray*}
\det(c_{ij})_{i,j=1}^r&=& \prod_{{i,j=1 \atop i<j
}}^{r}(a_i-a_j)(b_j-b_i)\prod_{i,j=1}^{r}\frac{1}{a_i-b_j}\\&&
-\frac{1}{a_r-b_{r}}\prod_{{i,j=1 \atop i<j
}}^{r-1}(a_i-a_j)(b_j-b_i)\prod_{i,j=1}^{r-1}\frac{1}{a_i-b_j}
\\ &=&
\frac{1}{a_r-b_{r}}\prod_{{i,j=1 \atop i<j
}}^{r-1}(a_i-a_j)(b_j-b_i)\prod_{i,j=1}^{r-1}\frac{1}{a_i-b_j}(
M-1).
\end{eqnarray*}
It follows that
$$\det(c_{ij})_{i,j=1}^r< 0.$$

\noindent\textbf{Class III.} {$c_{ij}\neq 0$ except for $c_{rr},\,
c_{r,r-1}\, \mbox{and}\, c_{r-1,r}$.}

In this case, we have $a_r>b_{r-2}$ but $a_r<b_{r-1}$,
$a_{r-2}>b_r$ but $a_{r-1}<b_r$. Recall that the \emph{rank} of a
matrix is the maximum number of linearly independent rows or
columns of the matrix. For a matrix $M=(m_{ij})_{i,j=1}^r$, let
$M^*$ be the matrix $(M_{ji})_{i,j=1}^r$, where $M_{ij}$ is the
cofactor of $m_{ij}$ in the expansion
$\det(M)=\sum_{i=1}^rm_{ij}M_{ij}$. Recall the following property:

\begin{equation}\label{bd-prop}
{\rm rank}(M^{*})=\left\{
\begin{array}{ll}
r, & \mbox{ if ${\rm rank}(M)=r$}\\[8pt]
1, & \mbox{ if ${\rm rank}(M)=r-1$}\\[8pt]
0, & \mbox{ if ${\rm rank}(M)<r-1$}
\end{array}
\right.
\end{equation}

We now consider the rank of $C^*=(C_{ji})_{i,j=1}^r$ where
$C_{ij}$ is the cofactor of $m_{ij}$ in the expansion
$\det(C(a,b))=\sum_{i=1}^rc_{ij}C_{ij}$. Recall that the minor
$C_{rr}$ is the determinant of the submatrix obtained from
$C(a,b)$ by deleting row $r$ and column $r$,  which turns out to
be  the restricted Cauchy matrix of Class \textbf{I}, and the
underlying matrices of $C_{r-1, r-1},\, C_{r, r-1},\,C_{r-1, r}$
are the restricted Cauchy matrices of Class \textbf{II}. Thus we
have
$$C_{r, r}>0,\quad C_{r-1, r-1}<0,\quad C_{r, r-1}>0 \quad \mbox{and} \quad C_{r-1, r}>0.$$
This implies that ${\rm rank}(C^*)\geq 2$. Hence ${\rm
rank}(C(a,b))=r$ because of \eqref{bd-prop}, namely
$\det(c_{ij})_{i,j=1}^r\neq 0$.

\noindent\textbf{Class IV.} {For all $i\leq r$, $j\leq r-1$ we
have $c_{ij}\neq 0$; $c_{1r}\neq 0$, and  $c_{2r}\neq 0$;
$c_{ir}=0$ if $i>2$.}

In this case, the restricted Cauchy matrix has the form
$$
(c_{ij})_{i,j=1}^r=
\begin{pmatrix}
\displaystyle\frac{1}{a_1-b_1} & \ldots & \displaystyle\frac{1}{a_1-b_{r-1}} & \displaystyle\frac{1}{a_1-b_r}\\
\displaystyle\frac{1}{a_2-b_1} & \ldots & \displaystyle\frac{1}{a_2-b_{r-1}} & \displaystyle\frac{1}{a_2-b_r}\\
\displaystyle\frac{1}{a_3-b_1} & \ldots & \displaystyle\frac{1}{a_3-b_{r-1}} & 0\\
\vdots            & \ldots & \vdots                & \vdots\\
\displaystyle\frac{1}{a_r-b_1} & \ldots &
\displaystyle\frac{1}{a_r-b_{r-1}} & 0
\end{pmatrix}.
$$
Expanding  along the last column, we get
\begin{eqnarray*}
\det(c_{ij})_{i,j=1}^r & = & (-1)^{r+1} \frac{1}{a_1-b_r}
\prod_{2\leq i<j\leq r} (a_i-a_j)\prod_{\atop {1\leq i<j\leq
r-1}}(b_j-b_i) \prod_{{2\leq i\leq r}\atop {1\leq j\leq r-1}}
\frac{1}{a_i-b_j}
\\&&+ (-1)^{r+2} \frac{1}{a_2-b_r} \prod_{{i\neq 2, j\neq 2}\atop {1\leq i<j\leq r}} (a_i-a_j)\prod_{\atop {1\leq i<j\leq r-1}}(b_j-b_i) \prod_{{i\neq 2}\atop
{1\leq j\leq r-1}} \frac{1}{a_i-b_j}
\\ & = & (-1)^{r+1}\prod_{{1\leq i<j\leq r}}(a_i-a_j)\prod_{\atop {1\leq i<j\leq r-1}}(b_j-b_i)\prod_{{i,j=1}\atop {j\neq r}}^{r}
\frac{1}{a_i-b_j}N,
\end{eqnarray*}
where
$$N=\frac{f(a_1)-f(a_2)}{a_1-a_2}$$
and $$
f(x)=\frac{(x-b_1)(x-b_2)\cdots(x-b_{r-1})}{(x-b_r)(x-a_3)\cdots(x-a_r)}.$$
Let $\delta=a_1-a_2$. We  obtain
\begin{eqnarray*}
\frac{f(a_1)}{f(a_2)} & = &
\frac{\displaystyle\frac{(a_1-b_1)(a_1-b_2)\cdots(a_1-b_{r-1})}{(a_1-b_r)(a_1-a_3)\cdots(a_1-a_r)}}
{\displaystyle\frac{(a_2-b_1)(a_2-b_2)\cdots(a_2-b_{r-1})}{(a_2-b_r)(a_2-a_3)\cdots(a_2-a_r)}}\\[5pt]
&=&
\frac{\displaystyle\frac{(a_1-b_1)}{(a_2-b_1)}\frac{(a_1-b_2)}{(a_2-b_2)}\cdots\frac{(a_1-b_{r-1})}{(a_2-b_{r-1})}}{\displaystyle
\frac{(a_1-b_r)}{(a_2-b_r)}\frac{(a_1-a_3)}{(a_2-a_3)}\cdots\frac{(a_1-a_{r})}{(a_2-a_{r})}}\\[5pt]
&=&
\frac{\displaystyle\frac{(\delta+a_2-b_1)}{(a_2-b_1)}\frac{(\delta+a_2-b_2)}{(a_2-b_2)}\cdots\frac{(\delta+a_2-b_{r-1})}{(a_2-b_{r-1})}}{
\displaystyle\frac{(\delta+a_2-b_r)}{(a_2-b_r)}\frac{(\delta+a_2-a_3)}{(a_2-a_3)}\cdots\frac{(\delta+a_2-a_{r})}{(a_2-a_{r})}}.
\end{eqnarray*}
Let $s\in\{b_1, \ldots, b_{r-1}\}$ and $s'\in\{a_3, \ldots, a_{r},
b_r\}$.  Then we have $s<s'$ and
\begin{equation}
\frac{(\delta+a_2-s)}{(a_2-s)}<\frac{(\delta+a_2-s')}{(a_2-s')}.
\end{equation}
It follows that $f(a_1)<f(a_2)$, namely $N<0$. Thus we have
$\det(c_{ij})_{i,j=1}^r> 0$ if $r$ is even and
$\det(c_{ij})_{i,j=1}^r< 0$ if $r$ is odd.

\vspace{.2cm} \noindent{\bf Acknowledgments.} This work was done
under the auspices of the 973 Project on Mathematical
Mechanization, the Ministry of Education, the Ministry of Science
and Technology, and the National Science Foundation of China. We
thank Professor Richard Stanley for bringing this problem to our
attention and for valuable comments. We also thank the referee for
the very pertinent comments and suggestions which helped to
significantly improve this paper.


\begin{thebibliography}{1}

\bibitem{B1} C. Bessenrodt, On hooks of Young diagrams, \textit{Ann.
Combin.} \textbf{2} (1998), 103-110.

\bibitem{B2} C. Bessenrodt, On hooks of skew Young diagrams and bars, \textit{Ann.
Combin.} \textbf{5} (2001), 37-49.

\bibitem{CYY} William Y. C. Chen, G.-G. Yan, and Arthur L. B.
Yang, Transformations of border strips and Schur function
determinants, \textit{J. Algebraic Combin.}, to appear.

\bibitem{D} V. Drinfeld, Hopf algebras and the Yang-Baxter
equation, \textit{Soviet Math. Dokl.} \textbf{32} (1985), 254-258.

\bibitem{GV1} I. Gessel and G. Viennot, Binomial determinants, paths, and hook length formulae,
\textit{Adv. Math.} \textbf{58} (1985), 300-321.

\bibitem{GV2} I. Gessel and G. Viennot, Determinants, paths, and plane partitions, preprint,
1989; available at {\tt http://www.cs.brandeis.edu/\~{}ira}.

\bibitem{HG} A. M. Hamel and I. P. Goulden, Planar decompositions
of tableaux and Schur function determinants, \textit{European J.
Combin.} \textbf{16}, 461-477.

\bibitem{Muir} T. Muir, A Treatise on the Theory of Determinants,
revised and enlarged by W. H. Metzler, Dover, New York, 1960.

\bibitem{NT} M. Nazarov and V. Tarasov, On irreducibility of tensor products of Yangian
modules associated with skew Young diagrams, \textit{Duke Math.
J.} \textbf{112} (2002), 343-378.

\bibitem{S1} R. P. Stanley, Enumerative Combinatorics, vol. 2,
Cambridge University Press, New York/Cambridge, 1999.

\bibitem{S2} R. P. Stanley, The rank and minimal border strip decompositions of
a skew partition, \textit{J. Combin. Theory Ser. A} \textbf{100}
(2002), 349-375.

\bibitem{stanleyprv} R. P. Stanley, private communication.
\end{thebibliography}
\end{document}